\documentstyle[12pt]{article}
\newcommand{\const}{\mathop{\rm const}\limits}

\newcommand{\supp}{\mathop{\rm supp}\limits}

\newcommand{\diam}{\mathop{\rm diam}\limits}

\newcommand{\Var}{\mathop{\rm Var}\limits}

\newcommand{\mod}{\mathop{\rm mod}\limits}

\begin{document}

\begin{center}

{\bf UNIFORM CENTRAL LIMIT THEOREM } \\

\vspace{3mm}

{\bf FOR MARTINGALES.} \par

\vspace{4mm}

 {\bf L.Sirota }  \\

\vspace{3mm}

\vspace{4mm}

  Department of Mathematics and computer science. \\
 Bar-Ilan University, 84105, Ramat Gan, Israel.\\

\vspace{3mm}

E-mail: \ sirota3@bezeqint.net\\

\vspace{7mm}

                    {\sc Abstract.}\\

 \end{center}

 \vspace{3mm}

 We study some sufficient conditions imposed on the sequence of martingale differences (m.d.) in the separable Banach spaces
of continuous functions defined on the metric compact set for the Central Limit Theorem  in this space. \par
 We taking into account the classical entropy terms, and use the theory of the so - called Grand Lebesgue Spaces
of random variables having power and exponential decreasing tail of distribution. \par

\vspace{4mm}

\par

\vspace{4mm}

 \par

\vspace{4mm}

\section{Notations. Statement of problem.}

\vspace{3mm}

 Let \\

{\bf A.} \ $ B $ be Separable Banach Space  with a norm $ || \ \cdot \ || B, $ briefly: $ B = (B, || \ \cdot \ ||B), $
 equipped with Borelian sigma-algebra; \\

\vspace{3mm}

{\bf B.} \ $ (\Omega, F,P) $ be probability triple with expectation $ E; $ \\

\vspace{3mm}

{\bf C.} \ $  F(i), \ i = 0,1,2, \ldots  $ be filtration, i.e. monotonically non-decreasing flow (sequence) of sigma-subfields
of source sigma-field $  F  $ such that $  F(0) $ is trivial sigma-algebra:  $  F(0) = \{ \emptyset, \Omega  \}; $ \\

\vspace{3mm}

{\bf D.} \ A sequence of {\it centered} martingale differences (m.d.) $  \xi(i), \ i = 1,2,\ldots,  $ i.e. Borelian distributed r.v. with
values in the space $  B  $  such that

$$
{\bf E} \xi(i)/F(k) = 0, \ k < i; \hspace{6mm}  {\bf E} \xi(i)/F(i) = \xi(i) \ (\mod {\bf P}). \eqno(1.0)
$$
 The last  equalities (1.0) may be explained as follows. For all non-random elements $ b^* $ from dual (conjugate) space $ B^* $

$$
{\bf E} \ b^*(\xi(i))/F(k) = 0, \ k < i; \hspace{6mm}  {\bf E} \ b^*(\xi(i))/F(i) = b^*(\xi(i)) \ (\mod {\bf P}). \eqno(1.0a)
$$

 Denote

$$
S(n) = \sum_{i=1}^n  \xi(i);
$$
then $ (S(n), F(n)) $ is really $ \ B  \ $  space valued mean zero martingale. \par

\vspace{3mm}

{\bf Definition 1.1.} We will say that the martingale $ (S(n), F(n)) $  or simple $ S(n) $ satisfies the Central Limit Theorem (CLT)
in the Banach space $  B, $ if the sequence of distributions of a random variables

$$
\eta(n) :=   n^{-1/2} S(n) = n^{-1/2} \sum_{i=1}^n  \xi(i)  \eqno(1.1)
$$
i.e. under classical norming  sequence $  1/\sqrt{n}, $
converges weakly as $ n \to \infty $ to the non-zero Gaussian distributed r.v. $ \eta(\infty) $ in the space $  B. $ \par

 Evidently, $ {\bf E} \eta(\infty) = 0. $ The covariation operator $ R_{\eta(\infty) }(\cdot)  $ of the limiting r.v. $ \eta(\infty)  $ has a form

$$
R_{\eta(\infty) }  (b^*) \stackrel{def}{=}  {\bf E} (b^*(\eta(\infty))^2
=\lim_{n \to \infty} n^{-1} \sum_{i=1}^n R_{\eta(i) }(b^*). \eqno(1.2)
$$

\vspace{3mm}

{\it  We will suppose that the Gaussian r.v. $ \eta(\infty) $  there exists, belongs to the space $  B \ (\mod({\bf P}))  $
and that for all $ b^* \in B^* $ the one-dimensional r.v. $ b^*(\eta(n)) $ converge in distribution as $ n \to \infty $
to one for the r.v. } $ b^*(\eta(\infty)). $ \par

\vspace{3mm}

{\bf  It remains to establish only the weak compactness, i.e. in  Yu.V.Prokhorov's sense \cite{Prokhorov1}  of distributions of the r.v. }
$  \eta(n) $  {\bf in the space $ B $ to  deduce the CLT for the martingale } $ S(n) $ {\bf in this space.} \par

\vspace{3mm}

 Note that the one-dimensional CLT for martingales is described in the famous book \cite{Hall1}; see also
 \cite{Bae1}, \cite{Pizier1}, \cite{Pizier2}, \cite{Walk1} etc. The CLT in Banach spaces for independent variables is
 considered in monographs \cite{Dudley1}, \cite{Ledoux1}, \cite{Ostrovsky1} and in many articles. \par

\vspace{3mm}

\section{ Uniform CLT for martingales in entropy terms. }

\vspace{3mm}

{\it This section may be considered as a simplification of the article  \cite{Bae1}. } \par

\vspace{3mm}

 Let $ (X = \{x \},d) $ be compact metric space relative some distance (or semi-distance) $  d = d(x_1, x_2), $
and let $ \{\xi_i \} =  \{ \xi_i(x) \}, \ x \in X $ be centered martingale differences relative the index $  i  $
random processes (r.p.)  (fields, r.f.);  the continuity with probability one of ones it follows from conditions of a next theorem. \par

Denote as ordinary

$$
\eta_n(x) := n^{-1/2} \sum_{i=1}^n \xi_i(x), \eqno(2.0)
$$

$$
\psi(p) := \sup_i \sup_{x \in X} |\xi_i(x)|_p = \sup_i \sup_{x \in X} \left[ {\bf E} |\xi_i(x)|^p   \right]^{1/p}. \eqno(2.1)
$$

\vspace{3mm}

 {\it  We suppose in what follows that the introduced function $ \psi = \psi(p) $  is finite at last for some value $ B, \ B = \const > 2; $
 may be $ B = \infty; $ } \hspace{4mm} then evidently it is finite for all the values $  p  $ from the set $  [1, B). $ \par

\vspace{3mm}

 The function $ \psi(\cdot) $ is called in the theory of Grand Lebesgue Spaces (GLS) as a {\it natural function } for the family of  the
random variables $  \{ \xi_i(x) \}, \ i = 1,2,\ldots; \ x \in X. $ \par

 We recall here briefly  the definition and some simple properties of the so-called Grand Lebesgue spaces;   more detail
investigation of these spaces see in \cite{Fiorenza3}, \cite{Iwaniec2}, \cite{Kozachenko1}, \cite{Liflyand1}, \cite{Ostrovsky1},
\cite{Ostrovsky2}; see also reference therein.\par

  Recently  appear the so-called Grand Lebesgue Spaces $ GLS = G(\psi) =G\psi =
 G(\psi; A,B), \ A,B = \const, A \ge 1, A < B \le \infty, $ spaces consisting
 on all  the random variables (measurable functions) $ f: \Omega \to R $ with finite norms

$$
   ||f||G(\psi) \stackrel{def}{=} \sup_{p \in (A,B)} \left[ |f|_p /\psi(p) \right]. \eqno(2.2)
$$

  Here $ \psi(\cdot) $ is some continuous positive on the {\it open} interval
$ (A,B) $ function such that

$$
     \inf_{p \in (A,B)} \psi(p) > 0, \ \psi(p) = \infty, \ p \notin (A,B).
$$

 We can accept in the sequel $  A = \const = 2, $ and following $  B > 2, $
taking into account the application in the theory of CLT in the space
of all numerical (real or complex) continuous functions $  C(X). $ \par

 We will denote
$$
 \supp (\psi) \stackrel{def}{=} (A,B) = \{p: \psi(p) < \infty, \}
$$

The set of all $ \psi $  functions with support $ \supp (\psi)= (A,B) $ will be
denoted by $ \Psi(A,B). $ \par
  This spaces are rearrangement invariant, see \cite{Bennet1}, and  are used, for example, in the theory of probability
 \cite{Kozachenko1},  \cite{Ostrovsky1}, \cite{Ostrovsky2}; theory of Partial Differential Equations \cite{Fiorenza3},
 \cite{Iwaniec2};  functional analysis \cite{Fiorenza3}, \cite{Iwaniec2},  \cite{Liflyand1},
 \cite{Ostrovsky2}; theory of Fourier series, theory of martingales, mathematical statistics, theory of approximation etc.\par

   The function $ \psi(\cdot) $ introduced in (2.1) generated the bounded {\it  natural distances } $ d_{i} = d_{i}(x_1, x_2) $
and $ \overline{d}    = \overline{d}(x_1, x_2), \ x_1,x_2 \in X  $  (more exactly, semi-distances) on the set $  X: $

$$
d_{i}(x_1, x_2) \stackrel{def}{=} ||\xi_i(x_1) - \xi_i(x_2)||G\psi, \eqno(2.3)
$$
so that

$$
d_{i}(x_1, x_2) \le 2; \hspace{5mm}  \forall i \ \Rightarrow ||\xi_i(x_1) - \xi_i(x_2)||G\psi \le d_{i}(x_1,x_2);
$$

$$
\overline{d}(x_1,x_2) \stackrel{def}{=} \sup_n \left[ \sqrt{ n^{-1} \sum_{i=1}^n d^2_i(x_1,x_2)} \right].\eqno(2.3a)
$$

\vspace{4mm}

 Let us introduce for any subset $ V, \ V \subset X $  and for arbitrary semi - distance on the set $ X $ \hspace{5mm}
 $ d = d(x_1, x_2), \ x_1, x_2 \in X $ the so-called {\it entropy } $ H(V, d, \epsilon) = H(V, \epsilon) $
 as a natural logarithm of a minimal quantity $ N(V,d, \epsilon) = N(V,\epsilon) = N $ of a balls $ S(V, t, \epsilon), \ t \in V: $

$$
S(V, t, \epsilon) \stackrel{def}{=} \{s, s \in V, \ d(s,t) \le \epsilon \},
$$
which cover the set $ V $ (covering numbers):

$$
N = N(V,d,\epsilon) = \min \{M: \exists \{t_i \}, i = 1,2,…, M, \ t_i \in V, \ V
\subset \cup_{i=1}^M S(V, t_i, \epsilon ) \}, \eqno(2.4)
$$
and we denote also $  D = D(d) = \diam(X, d) = \sup_{x_1,x_2 \in X} d(x_1, x_2), $

$$
H(V,d,\epsilon) = \log N; \ S(t_0,\epsilon) \stackrel{def}{=}
 S(X, t_0, \epsilon), \ H(d, \epsilon) \stackrel{def}{=} H(X,d,\epsilon). \eqno(2.4a)
$$

 It follows from Hausdorff's theorem that
$ \forall \epsilon > 0 \ \Rightarrow H(V,d,\epsilon)< \infty $ iff the
metric space $ (V, d) $ is precompact set, i.e. is the bounded set with
compact closure.\par

\vspace{4mm}

 We will distinguish in the sequel  two cases: $  B < \infty,  $ {\it finite case}, and $  B = \infty, $ {\it  infinite case}. \par

\vspace{3mm}

{\bf Finite case.} \\

\vspace{3mm}

 The probabilistic Grand Lebesgue Spaces $ G\Psi(2,B)  $ with $ 2 < B < \infty  $ are in detail investigated in articles
 \cite{Liflyand1},  \cite{Ostrovsky108}, including  consideration many examples.\par

 We will use in this case the following inequality for the arbitrary sequence $  \{ \zeta_k \}  $ of martingale difference, see   the
famous article of A.Osekovski \cite{Osekowski1}, see also \cite{Ostrovsky604}:

$$
 \left|\ n^{-1/2} \ \sum_{k=1}^n \zeta_k \ \right|_p \le K_{Os} \cdot \frac{p}{\ln p}  \cdot \sqrt{ \sum_{k=1}^n |\zeta_k|_p^2/n  },
 \eqno(2.5)
$$
where  the "Osekowski's"  constant  $  K_{Os} $  is less than 15.5879. \par
 It is interest to note that  at the same estimate was before obtained by H.Rosenthal \cite{Rosenthal1} for independent variables;
in this case the exact value of this constant ("Rosenthal's constant")  $ C_R \approx 0.6535, $  see \cite{Ostrovsky601}. \par

\newpage

{\bf Infinite case.} \\

\vspace{3mm}

  The spaces  $  G\Psi(2, \infty)  $  are convenient, e.g., for the investigation of the random variables and vectors with
exponential decreasing tail of distribution.  Indeed,
if for some non-zero r.v. $ \xi \ $ we have $ 0 < ||\xi||G(\psi) < \infty, $
 then for all positive values $ u $

$$
{\bf P}(|\xi| > u) \le 2 \ \exp \left( - \overline{\psi}^*(\log x /||\xi||G(\psi))  \right),
\eqno(2.6)
$$
where  $ \overline{\psi}(p) = p \ \log \psi(p) $ and the symbol $ g^* $ denotes some modification
of the Young-Fenchel, or Legendre  transform of the function $ g: $

$$
g^*(y) = \sup_{x \ge 2} (x y - g(x)).
$$
 see \cite{Kozachenko1}, \cite{Ostrovsky1}, chapters 1,2.\par
 As a  consequence: if

 $$
\forall x > e^2 \ \Rightarrow  \overline{\psi}^*(\log x)  > 0,
 $$
then the space $ G\psi $ coincides with {\it exponential} Orlicz's  space over  our probabilistic space
$ (\Omega,F,{\bf P} ) $ with $ N- $ function of a form

$$
N(u) = \exp( \overline{\psi}^*(\log |u|) ), \ |u| > e^2; \ N(u) = C \cdot u^2, |u| \le e^2. \eqno(2.7)
$$

  Conversely: if a r.v. $ \xi $ satisfies (2.6), then $ \xi \in G\psi, \  ||\xi||G(\psi) < \infty. $ \par

\vspace{3mm}

{\bf Example 2.1.}  Introduce as a particular case the following norm for the r.v. $  \xi: $

$$
K :=  ||\xi||_{(q)} := \sup_{p \ge 2} \left[ \frac{|\xi|_p}{p^{1/q}} \right],   \ q = \const > 0;
$$
then $ K = ||\xi||_{(q)} \in (0,\infty) \Leftrightarrow $

$$
T(\xi,x) :=   \max ({\bf P}(\xi > x), {\bf P}(\xi < x)  ) \le \exp \left( - C(q) (x/K)^q  \right), \ x > 1. \eqno(2.8)
$$

\vspace{3mm}

 So, the theory of $  G\psi $ spaces of random variables gives a very convenient  apparatus
for investigation of a random variables with exponential decreasing tails of distribution. \par

\vspace{3mm}

 Let $ \psi \in \Psi(2,B);  $ we introduce the so-called Rosenthal's transform $  \psi_R(\cdot) $ as follows:

$$
\psi_R(p) := \frac{p}{\log p} \cdot \psi(p). \eqno(2.9)
$$

 Evidently,  $ \psi(p) \le e \cdot \psi(p);  $ and in addition  if $  B < \infty, $ then

$$
\psi_R(p) \asymp \psi(p),
$$
and this is not true if $  B = \infty. $\par

\vspace{3mm}

 Let us denote for arbitrary function $ \psi \in \Psi $

$$
\psi_*(x) := \inf_{y \in (0,1)} (x y + \log \psi(1/y)),
$$
and introduce the following functional ("entropy integral, covering integral")

$$
J(\psi, d) \stackrel{def}{=} \int_0^D \exp ( \psi_*( \log 2 + H(X,d,\epsilon)  ) ) \ d \epsilon. \eqno(2.10)
$$

\vspace{4mm}

{\bf Remark 2.1.} If we introduce the {\it discontinuous} function

$$
\psi_r(p) = 1, \ p = r; \psi_r(p) = \infty, \ p \ne r, \ p,r \in (A,B)
$$
and define formally  $ C/\infty = 0, \ C = \const \in R^1, $ then  the norm
in the space $ G(\psi_r) $ coincides with the $ L_r $ norm:

$$
||f||G(\psi_r) = |f|_r.
$$

 Thus, the Bilateral Grand Lebesgue spaces are direct generalization of the
classical exponential Orlicz's spaces  as well as of the classical Lebesgue-Riesz spaces $ L_r. $ \par

\vspace{3mm}

 We need to introduce again some notations. Let  $ T(x), \ x > 0 $ be a  {\it tail - function,}  i.e.  such that
$ T(0) = 1, \ T(\cdot) $ is monotonically
decreasing, right continuous and such that $ T(\infty) = 0. $ \par

 We denote for the tail-function $ T(\cdot) $  the following operator (non-linear)
$$
W[T](x) = \min \left(1, \inf_{v > 0} \left[ \exp(-x^2/(8v^2)) - \int_v^{\infty}
x^2 \ dT(x) \right] \right), \eqno(2.11)
$$
if there exists the second moment
$$
\int_0^{\infty} x^2 \ |dT(x)| < \infty.
$$

\vspace{3mm}

{\bf Lemma 2.1. } (See \cite{Ostrovsky603}).  Let $ \ \xi(i)  $  be a sequence of {\it centered} martingale-differences relative
to some filtration $ \{F(i)\} $ and $ T(\xi(i),x) \le T(x), \ T(x) $ be some tail-function. Then at $ x > 1 $

$$
\sup_{b: \sum b^2(i) =1} T \left(\sum_i b(i) \xi(i), x \right)  \le W[T](x). \eqno(2.12)
$$
 In particular,

$$
 \sup_n T \left( n^{-1/2} \sum_{i=1}^n  \xi(i), x \right)  \le W[T](x). \eqno(2.12a)
$$

 If for instance

 $$
 \exists q,K = \const > 0 \ \Rightarrow   T(x) \le \exp \left(-(x/K)^q \right), \ x \ge 0,
 $$
then

$$
 \sup_n T \left( n^{-1/2} \sum_{i=1}^n \xi(i), x \right) \le \exp \left( - C(q) (x/K)^{ 2q/(2 + q) } \right). \eqno(2.12b)
$$

 Denote also

$$
\sigma^2  = \sigma^2( \{  \xi_i(\cdot) \}  ) = \inf_{x \in X} \sup_n \left[ n^{-1} \sum_{k=1}^n \Var(\xi_k(x)) \right]. \eqno(2.13)
$$

\vspace{4mm}

{\bf Theorem 2.1.} ("Power" level.) \par
Suppose that for our sequence of functional  martingale  differences $ \{  \xi_i(\cdot) \}  $

$$
 \sigma^2( \{  \xi_i(\cdot) \}  )  < \infty \eqno(2.14)
$$
and

$$
J(\psi_R, \overline{d} ) = \int_0^{D(\overline{d})} \exp \left( \psi_{R,*} (\ln 2 + H(X, \overline{d}, \epsilon))  \right)
 \ d \epsilon < \infty. \eqno(2.15)
$$

 Then the family of distributions of the sequence of random fields $ \{ \eta_n(\cdot) \} $ is weakly compact in the space
 $  C(X, \overline{d}). $\par

 \vspace{3mm}

{\bf Proof.}\\
{\bf 0.} As long as $  \psi(p) \le  e\cdot \psi_R(p), $  and

$$
d_k(x_1, x_2) \le C \overline{d}(x_1, x_2),
$$
it follows from the condition (2.15) the convergence of the following
entropy integral for the individual  r.f.

$$
J(\psi, \overline{d} ) = \int_0^{D(\overline{d})} \exp \left( \psi_{*} (\ln 2 + H(X, \overline{d}, \epsilon))  \right)
 \ d \epsilon < \infty. \eqno(2.15a)
$$
 Therefore, each r.f. $ \xi_k(x) $ is continuous with probability one relative the distance $  \overline{d}(x_1, x_2). $\\

{\bf 1.}  The condition (2.14) imply that there exists at least  one point $ x_0 $ in the set $ X $ for which

$$
\sup_n \left[ n^{-1} \sum_{k=1}^n \Var(\xi_k(x_0)) \right] < \infty.
$$

 Since the martingale differences are centered and non-correlated,

$$
\sup_n {\bf E} \left[ \eta_n(x_0) \right]^2 < \infty.
$$
 Therefore, the family of distributions on real line of one-dimensional r.v.  $ \{ \eta_n(x_0)  \} $  is weakly compact.\par

\vspace{3mm}

{\bf 2.} We apply the inequality (2.5) for the variables $ \zeta_k = \xi_k(x_1) - \xi_k(x_2),  $ where $  x_1, x_2  $ are non-random
elements of the set $ X $ such that $ \overline{d}(x_1, x_2) > 0 $  (the case $ \overline{d}(x_1, x_2) = 0 $ is trivial):

$$
|\eta_n(x_1) - \eta_n(x_2)|_p \le K_{Os} \cdot \frac{p}{\ln p } \cdot \sqrt{ \left[ n^{-1} \sum_{k=1}^n |\xi_k(x_1) - \xi_k(x_2)|_p^2  \right] } \le
$$

$$
K_{Os} \cdot \frac{p}{\ln p }  \cdot \psi(p) \cdot \overline{d}(x_1, x_2) = K_{Os} \cdot \psi_R(p) \cdot \overline{d}(x_1, x_2), \eqno(2.16)
$$
or equally

$$
\sup_n || \eta_n(x_1) - \eta_n(x_2)||G\psi_R \le  K_{Os} \cdot \overline{d}(x_1, x_2). \eqno(2.16a)
$$

  The statement of theorem 2.1 it follows from theorem 4.4.2 of the monograph \cite{Ostrovsky1}, chapter 4, section 4. \par

\vspace{3mm}

 {\bf Example 2.1.} Let $ \psi(p) = \psi_r(p) $ for some $ r = \const \ge 2; $ then the space $ G\psi_r $ coincides with
the classical Lebesgue - Riesz space $ L_r = L_r(\Omega, {\bf P}). $ \par
 Then the condition (2.15) is equal to the following famous G.Pizier's condition

$$
\int_0^{D(d_r)} N^{1/r}(X, d_r,\epsilon) \ d \epsilon < \infty,
$$
i.e. as in the independent case, see e.g. \cite{Pizier2}. Here $ d_r(\cdot, \cdot) $ is the so-called natural Pizier
distance

$$
d_r(x_1, x_2) = \sup_k |  \xi_k(x_1) - \xi_k(x_2)|_r.
$$

\vspace{3mm}

 Let us consider now the so-called "exponential level". Let again $ q = \const > 0; $ define the following distance  on the set $  X  $

$$
\rho_q(x_1, x_2) := \sup_k || \xi_k(x_1) - \xi_k(x_2)||_{(q)}. \eqno(2.17)
$$

 \vspace{3mm}

 {\it It will be presumed the boundedness of this distance, as well as the finiteness of the variable }
$  \sigma^2( \{  \xi_i(\cdot) \}  ) . $ \par

\vspace{3mm}

{\bf Theorem 2.2.} Suppose in addition

$$
\int_0^{D(\rho_q)} H^{ (2 + q)/(2 q)  }(X, \rho_q, \epsilon) \ d \epsilon < \infty. \eqno(2.18)
$$

Then the family of distributions of the sequence of random fields $ \{ \eta_n(\cdot) \} $ is weakly compact in the space
 $  C(X, \rho_q). $\par

 \vspace{3mm}

 {\bf Proof.} We apply the inequality (2.18) and described below the properties of Grand Lebesgue Spaces:

 $$
 \sup_n || \ \eta_n(x_1) - \eta_n(x_2) \ ||_{( 2q/(q + 2) )} \le C_1(q) \rho_q(x_1, x_2). \eqno(2.19)
 $$
 The proposition of theorem 2.2 follows immediately from theorem 4.3.2 of the book \cite{Ostrovsky1}, chapter 4, section 4.3;
 also \cite{Kozachenko1}. \par

 \vspace{4mm}

 {\bf Remark 2.1.} In the {\it  independent case, } i.e. when the r.f. $ \xi_k(x) $ are in addition commonly independent,
 the condition (2.18) looks as follows:

$$
\int_0^{D(\rho_q)} H^{ 1/ q  }(X, \rho_q, \epsilon) \ d \epsilon < \infty,
$$
see \cite{Kozachenko1}; i.e. is unlike in comparison  to the general (martingale) case, in contradiction to the "power case".\par
 The eventually value of power of the variable  $ H(X, \rho_q, \epsilon) $ in the condition (2.18)  for martingales is now unknown. \par

\vspace{3mm}

{\bf Example 2.2.} The condition (2.18) is satisfied if for example $  X $ is bounded closed subset of the Euclidean space $  R^d, \ d = 1,2,\ldots $
and

$$
\rho_{q}(x_1, x_2) \le C_1 \ |x_1 - x_2|^{\alpha}, \ x_1,x_2 \in X, \  0 < \alpha = \const  \le 1.
$$
 In this case

$$
N(X, \rho_q, \epsilon) \le C_2(\alpha,q) \ \epsilon^{-d/\alpha}, \ \epsilon > 0.
$$

\vspace{4mm}

\section{Concluding remarks.}

\vspace{3mm}

 It is interesting by our opinion to obtain the CLT for functional martingales
in the space $  C(X,d) $ in more modern terms of majorizing measures
in the spirit of the article of B.Heinkel \cite{Heinkel1}, see also \cite{Fernique1}, \cite{Fernique2}, \cite{Fernique3},
\cite{Talagrand4}, \cite{Talagrand5}; as well as in the H\"older's space, see \cite{Klicnarov'a1},
\cite{Ratchkauskas1}, \cite{Ratchkauskas2}, by means of martingale inequalities \cite{Barlow1}, \cite{Ostrovsky603}.\par

\end{document}